\documentclass{amsart}

\usepackage[T1]{fontenc}
\usepackage[british]{babel}
\usepackage{amssymb}

\usepackage{mfmath}

\numberwithin{equation}{section}
\numberwithin{theorem}{section}
\allowdisplaybreaks[1]

\def\and{\qquad\text{and}\qquad}

\def\Ll{\boldsymbol\Lambda}
\def\Sl{\mathbf S}

\def\Su{\Sl^*}
\def\Lu{\Ll^{\!*}}
\def\Llprime{\Ll^{\!\prime}}

\def\LMod{\hbox{\bf\em M\kern-0.1em od}}
\def\lMod#1{#1\hbox{-\bf\em M\kern-0.05em od}}
\def\rComod#1{\hbox{{\bf\em Comod}-}#1}
\def\wMod#1{\lMod{B\Sl}}
\def\wComod#1{\rComod{B\Ll}}

\def\U{\mathcal U}

\def\otimesunder #1{\mathbin{\mathchoice
  {\mathop\otimes\limits_{\mkern-20mu #1\mkern-20mu}}%
  {\otimes_{#1}}{\otimes_{#1}}{\otimes_{#1}}}}

\def\tildeotimes{\mathbin{\tilde\otimes}}

\def\facet{\mathrel{<_1}}
\DeclareMathOperator{\Or}{or}
\DeclareMathOperator{\lin}{lin}
\DeclareMathOperator{\pr}{pr}

\def\Ring{\Z}
\def\SR#1{\Ring[#1]}
\def\SRco#1{\Ring(#1)}

\def\Hcld{H^{\rm cld}}

\title[On the integral cohomology of smooth toric varieties]%
  {On the integral cohomology \\ of smooth toric varieties}
\author{Matthias Franz}
\address{Institut Fourier, Universit\'e de Grenoble~I, BP~74, 38402~Saint-Martin d'H\`eres, France}
\email{matthias.franz@ujf-grenoble.fr}


\begin{document}

\begin{abstract}
  Let $X_\Sigma$ be a smooth, not necessarily compact toric variety.
  We show that a certain complex, defined in terms of the fan~$\Sigma$,
  computes the integral cohomology of~$X_\Sigma$,
  including the module structure over the homology of the torus.
  In some cases we can also give the product.
  As a corollary we obtain that the cycle map from Chow groups
  to integral Borel--Moore homology is split injective
  for smooth toric varieties.
  Another result is that the differential algebra of singular cochains
  on the Borel construction of~$X_\Sigma$ is formal.
\end{abstract}

\maketitle

\section{Introduction}

Let $X=X_\Sigma$ be a smooth toric variety, associated with a regular
fan~$\Sigma$ in the rational vector space spanned by a lattice~$N\cong\Z^r$.
The algebraic torus~$T=T_N\cong(\C^*)^r$ acts on~$X$.
The equivariant cohomology~$H_T^*(X)$ of~$X$ is by definition
the cohomology of its Borel construction~$X_T=(ET\times X)/T$.
(Here, as in the rest of the paper, we use integer coefficients.)
It is known that $H_T^*(X)$ 
is isomorphic to the Stanley--Reisner ring~$\SR\Sigma$, 
which only depends on the combinatorics of~$\Sigma$
(Bifet--DeConcini--Procesi~\cite{BifetDeConciniProcesi:90},
Davis--Januszkiewicz~\cite{DavisJanuszkiewicz:91}, Brion~\cite{Brion:96}).
(See Section~\ref{preliminaries} for a definition of~$\SR\Sigma$.)

We begin with a strengthening of this result.
We denote by~$C^*$ the normalised singular cochain functor
(with integer coefficients).
Recall that a differential algebra is called formal if it is quasi-isomorphic
to its homology.

\begin{theorem}\label{Stanley-Reisner-formal}
  The differential algebra~$C^*(X_T)$ is
  quasi-isomorphic to~$\SR\Sigma$.
  In particular, it is formal.
\end{theorem}

The proof of Theorem~\ref{Stanley-Reisner-formal} also leads to 
a description of the integral cohomology of~$X$.
Fix a basis~$(x_i)$ of~$N$ with dual basis~$(\xi_i)$.
The homology~$H(T)=\Ll$ is an exterior algebra, generated by the~$x_i$.
Recall that the $T$-action on~$X$ induces an action of~$\Ll$ on~$H^*(X)$.
For~$X=T$, this action is the contraction of forms in~$\Lu=H^*(T)$
by elements in~$\Ll$.
Moreover, the cohomology~$\Su=H^*(BT)$ of the classifying space of~$T$
is a symmetric algebra, generated by the~$\xi_i$,
and $\SR\Sigma$ is canonically an $\Su$-module.
(This action is not determined by the combinatorics of~$\Sigma$, though.)

\begin{theorem}\label{cohomology-module}
  The cohomology~$H^*(X)$ of~$X$ is isomorphic to
  the homology of the differential $\Ll$-module
  \begin{equation}\label{complex-cohomology}
    \Lu\otimes\SR\Sigma,
    \qquad
    d(\alpha\otimes\sigma)=\sum_{i=1}^r x_i\cdot\alpha\otimes\xi_i\sigma.
  \end{equation}
\end{theorem}

On~$\Lu\otimes\SR\Sigma$ there is the canonical product
\begin{equation}
  \label{canonical-product}
  (\alpha\otimes\sigma)\cdot(\alpha'\otimes\sigma')
    =\alpha\wedge\alpha'\otimes\sigma\sigma'
\end{equation}
(which is well-defined because the generators of~$\Ll$
act as derivations on~$\Lu$).
I conjecture that it always gives the correct product in cohomology,
but I can only prove it in a special case:

\begin{theorem}\label{cohomology-algebra}
  Suppose
  that $X$ is a toric subvariety of the projective variety~$(\mathbf P^1)^r$.
  Then the product~\eqref{canonical-product}
  induces the cup product in cohomology.
\end{theorem}

These results are related to work
of Buchstaber and Panov~\cite[Ch.~7]{BuchstaberPanov:02},
but their methods do not work over the integers~\cite{Panov}.
Theorem~\ref{cohomology-algebra} applies in particular
to complements of complex coordinate subspace arrangements.
We therefore answer a question of
Buchstaber--Panov~\cite[Problem~8.14]{BuchstaberPanov:02}.

Since $X$ is smooth, its cohomology~$H^*(X)$ is Poincar\'e dual
to its homology with closed support~$\Hcld(X)$ (Borel--Moore homology).
To describe a complex dual to~\eqref{complex-cohomology},
we fix an orientation of each cone in~$\Sigma$.
For~$\sigma$,~$\tau\in\Sigma$, $\sigma$ a facet of~$\tau$
(written: $\sigma\facet\tau$),
set~$\Or_{\sigma\tau}=\pm1$, depending on whether
the chosen orientations of $\sigma$~and~$\tau$ are compatible or not,
as in simplicial homology.
(See Section~\ref{proof-Borel-Moore-homology} for details on this
as well as on the cap product used below.)
Denote by~$\Ll^\sigma$ the exterior
algebra on the quotient of~$N$ by the sublattice
spanned by the rays of~$\sigma$,
and by~$\pr_\sigma$
the projection~$\Ll^\tau\to\Ll^\sigma$ for a face~$\tau$ of~$\sigma$,
and also the projection~$\Ll\to\Ll^\sigma$.

\begin{proposition}\label{Borel-Moore-homology}
  The double complex~$C(\Sigma)$ with
  \begin{subequations}
  \begin{equation}
    C_{pq}(\Sigma)=\bigoplus_{\sigma\in\Sigma_{r-p}}\Ll^\sigma_q,
  \end{equation}
  and
  \begin{equation}
    d(a)=\sum_{\sigma\facet\tau}\Or_{\sigma\tau}(-1)^{\degree a}\pr_\tau(a),
    \qquad
    x_i\cdot a=\pr_\sigma(x_i)\wedge a
  \end{equation}
  \end{subequations}
  for~$a\in\Ll^\sigma$, is quasi-isomorphic
  to the dual of~\eqref{complex-cohomology} as differential $\Ll$-module.
  In particular, it computes the Borel--Moore homology~$\Hcld(X)$ of~$X$.
\end{proposition}

This complex arises as the $E^1$~term of the spectral sequence
for~$\Hcld(X)$ coming from the filtration of~$X$ by orbit dimension,
\cf~\cite{Totaro:??},~\cite{Fischli:92},~\cite{Jordan:97}.
From this point of view, Proposition~\ref{Borel-Moore-homology}
says that this spectral sequence degenerates on the $E^2$~level
and moreover that there is no extension problem.

Since the diagonal terms~$H_{pp}(\Sigma)$ of the homology of~$C(\Sigma)$
are known to be the Chow groups of~$X$
(Fulton--Sturmfels~\cite{FultonSturmfels:97}),
we can partially answer a question of Totaro~\cite[Sec.~8]{Totaro:??}:

\begin{corollary}\label{Chow-injective}
  For smooth toric varieties, the cycle map
  from Chow groups to integral Borel--Moore homology is split injective.
\end{corollary}

We have stated all results for integer coefficients,
but in fact they hold for an arbitrary commutative ring with unit.
For rational or real coefficients, most of our results are already known.
Moreover, the proofs would much easier in this case,
for example because one can use (commutative) differential forms
instead of (non-commutative) singular cochains.

Our approach is based (or at least inspired) by the relation
between Koszul duality and equivariant cohomology
as described in~\cite{Franz:03}.
Theorems \ref{Stanley-Reisner-formal}~and~%
\ref{cohomology-module} appeared,
with slightly different proofs, in the author's thesis~\cite{Franz:01}.

\begin{acknowledgements}
  The author thanks Andrzej Weber for helpful discussions.
\end{acknowledgements}

\section{Preliminaries}
\label{preliminaries}

We shall mostly work in the homological setting.
The results stated in the introduction then follow by dualisation.
Most of the time we could equally use cohomology, but for our proof
of Theorem~\ref{cohomology-module} we need the homological results.

\def\Ring{R}
We work over a commutative ring~$\Ring$ with unit.
All (co)algebras and (co)modules will be differential graded;
differentials always lower degrees.
In particular, cohomology is negatively graded.

Let $N\simeq\Z^r$ be a lattice.
We write $\Ll$ for the exterior algebra over~$N$.
The symmetric coalgebra over~$N$ is denoted by~$\Sl$
and its dual, the symmetric algebra over~$N^*$, by~$\Su$.

We call a coalgebra together with a coalgebra map to~$\Sl$
an $\Sl$-coalgebra. A $\Ll$-module that is also a coalgebra
with a $\Ll$-equivariant diagonal is called a $\Ll$-coalgebra.
(Recall that $\Ll$ is a Hopf algebra
with all elements in~$\Ll_1=N\otimes\Ring$ being primitive.)
A morphism of $\Sl$-coalgebras or $\Ll$-coalgebras is a morphism
of coalgebras compatible with the additional structure.

Any choice of basis~$(x_i)$ of~$N$
gives algebra generators (of degree~$1$) of~$\Ll$
and coalgebra cogenerators (of degree~$2$) of~$\Sl$.
We write~$x_i$ for both of them, and $\xi_i$ for the dual basis elements,
considered as generators of~$\Lu$~or~$\Su$.
In terms of theses bases, the action of~$\Ll$ on~$\Lu$ is given by
$$
  x_i\cdot(\xi_{j_1}\wedge\cdots\wedge\xi_{j_q})
  =(-1)^{k-1}\,
    \xi_{j_1}\wedge\cdots\wedge\widehat\xi_{j_k}\wedge\cdots\wedge\xi_{j_q}
$$
if $j_k=i$ (and all other indices different from~$i$).
Note that any $\Sl$-comodule~$M$ is also an $\Su$-module.
We write~$\sigma\cap m$ for the induced action
of~$\sigma\in\Su$ on~$m\in M$.

For any $\Ll$-module~$L$ let $\Sl\tildeotimes L$
be the complex~$\Sl\otimes L$ with differential
$$
  d(s\otimes l)=s\otimes d l-\sum_{i=1}^r\xi_i\cap s\otimes x_i\cdot l.
$$
Note that the differential is independent of the choses basis~$(x_i)$.
(This is the Koszul dual of~$N$, \cf~\cite[Sec.~2.6]{Franz:03}.
We have introduced the minus sign to end up with the same
differential~\eqref{complex-cohomology} as in~\cite{BuchstaberPanov:02}.)

\smallbreak

The normalised singular chain and cochain functors
with coefficients in~$\Ring$
are denoted $C$~and~$C^*$, respectively.

\smallbreak

A cone is always rational and also pointed,
\ie, does not contain any line.
By a ray, we mean a one-dimensional cone. The rays of a cone
are its one-dimensional faces.
We write $x_\rho$ for the shortest vector with integral coordinates
contained in the ray~$\rho$.

Throughout this paper, $\Sigma$ denotes a fan in~$N$,
which is assumed to be regular,
\ie, for any cone~$\sigma\in\Sigma$,
the set~$\{x_\rho:\hbox{$\rho$ a ray of~$\sigma$}\}$
can be extended to a $\Z$-basis of~$N$.
We write~$\Sigma_p$ for the $p$-dimensional cones in~$\Sigma$.

The Stanley--Reisner ring~$\SR\Sigma$ is the polynomial algebra 
on generators~$\xi_\rho$, $\rho\in\Sigma_1$,
divided by all monomials~$\xi_{\rho_1}\cdots\xi_{\rho_k}$
such that $\{\xi_{\rho_1},\ldots\xi_{\rho_k}\}$ are not
the rays of a cone in~$\Sigma$.

Dually, the Stanley--Reisner coalgebra~$\SRco\Sigma$ is the coalgebra
with $\Ring$-basis, all monomials~$\rho_1\cdots\rho_k$ such that
$\rho_1$,~\ldots,~$\rho_k$ are the rays of some~$\sigma\in\Sigma$.
The~$\rho_i$ need not be distinct, but their order is irrelevant.
Each~$\rho$ has degree~$2$, and the diagonal decomposes each monomial
in all possible ways.

The map that sends each~$\rho\in\SRco\Sigma$ to~$x_\rho\in N$
extends to a morphism of coalgebras~$\SRco\Sigma\to\Sl$.
The dual algebra morphism $\Su\to\SR\Sigma$ is given by
$$
  \xi_i\mapsto\sum_{\rho\in\Sigma}\xi_i(x_\rho)\,\xi_\rho.
$$

The standard reference for toric varieties is~\cite{Fulton:93}.
We write $X_\sigma\subset X_\Sigma$
for the open affine subvariety determined by~$\sigma\in\Sigma$.

The choice of basis for~$N$ also determines 
a decomposition~$T\cong(\C^*)^r$.
The map that sends each~$x_i$ to the loop in~$T$ around the $i$-th factor
extends to a quasi-isomorphism of Hopf algebras
$$
  \Ll\to C(T),
$$
\cf~\cite[Lemma~4.1]{Franz:03}.
In particular, each $C(X_\sigma)$, $\sigma\in\Sigma$
becomes a $\Ll$-module this way
such that the induced $\Ll$-action on~$H(X_\sigma)$ is the standard one.

\section{Proof of Theorem~\protect\ref{Stanley-Reisner-formal}}\label{toric}

\begin{proposition}\label{subfan-positive-orthant-formal}
  Let $X$ be a toric subvariety of~$(\mathbf P^1)^r$.
  Then $\Sl\tildeotimes C(X)$ and $\SRco\Sigma$ are quasi-isomorphic
  as $\Sl$-coalgebras.
\end{proposition}

Note that the assumption on~$X$ means that each cone in~$\Sigma$ is generated
by a subset of~$\{\pm x_i\}$ where $(x_i)$ is a basis of the lattice~$N$.

\begin{proof}
  We will construct a sequence of quasi-isomorphisms
  leading from~$\Sl\otimes C(X)$ to~$\SRco\Sigma$.
  Cover $X$ by the toric subvarieties~$X_{\sigma_i}$
  corresponding to the maximal cones~$\sigma_1$,~\ldots,~$\sigma_m\in\Sigma$
  and call this covering $\U$.
  In a first step, we replace $\Sl\tildeotimes C(X)$
  by~$\Sl\tildeotimes E(\U)$.

  Here $E(\U)$ denotes the Mayer--Vietoris double complex
  associated with~$\U$. (See for example~\cite[\S 15]{BottTu:82}
  for the dual construction with cochains,
  or~\cite[\S A.2--A.6]{GugenheimMay:74}.)
  In order to describe $E(\U)$, it is useful to introduce some abbreviations.
  We will write $C(j_0,\ldots,j_p)$
  for~$C(X_{\sigma_{j_0}}\cap\cdots\cap X_{\sigma_{j_p}})%
  =C(X_{\sigma_{j_0}\cap\cdots\cap\sigma_{j_p}})$
  and $\iota(j_0,\ldots,j_p)$ for any inclusion
  into~$X_{\sigma_{j_0}\cap\cdots\cap\sigma_{j_p}}$.
  The double complex~$E(\U)$ has components
  $$
    E_{pq}(\U)
    =\bigoplus_{1\le j_0<\cdots<j_p\le m}
       C_q(j_0,\ldots,j_p)
  $$
  and differential
  $$
    d_{E(\U)}(c)
    =(-1)^p d_{C(j_0,\ldots,j_p)}(c)
    +\sum_{k=0}^p(-1)^k\iota(j_0,\ldots,\widehat{j_k},\ldots,j_p)_*(c)
  $$
  for~$c\in C(j_0,\ldots,j_p)$.
  (Interpret $\iota()_*(c)$ as~$0$.)
  Moreover, it is a $\Ll$-coalgebra with $\Ll$-action
  $$
    x_i * c=(-1)^p\,x_i\cdot c,
  $$
  where the action on the right is that
  of~$C_q(j_0,\ldots,j_p)$,
  and diagonal
  $$
    \Delta_{E(\U)}(c)=\sum_{k=0}^p\pm
      \bigl(\iota(j_0,\ldots,j_k)_*\otimes\iota(j_k,\ldots,j_p)_*\bigr)
      \Delta_{C(j_0,\ldots,j_p)}(c),
  $$
  where the sign is~$(-1)^{(p-k)\degree{c'}}$ for each term~$c'\otimes c''$
  appearing in~$\Delta_{C(j_0,\ldots,j_p)}(c)$.
  The canonical map~$E(\U)\to C(X)$ given
  by the inclusion~$C(j_0)\hookrightarrow C(X)$ for~$p=0$
  (and the zero map otherwise)
  is a quasi-isomorphism of $\Ll$-coalgebras.
  Therefore, the induced map~$\Sl\tildeotimes E(\U)\to\Sl\tildeotimes C(X)$
  is a quasi-isomorphism of coalgebras.
  In addition, it endows $\Sl\tildeotimes E(\U)$ with a coalgebra map to~$\Sl$.

  More important than the details of the construction of~$E(\U)$ is the fact
  that the complex~$\Sl\tildeotimes E(\U)$
  is made up of $\Sl$-coalgebras~$\Sl\tildeotimes C(X_\sigma)$,
  connected by maps induced by inclusions~$X_\tau\hookrightarrow X_\sigma$,
  and that this construction is natural with respect
  to maps of $\Sl$-coalgebras which are compatible with the maps giving
  the differentials. Moreover, a standard spectral sequence argument
  shows that replacing each coalgebra by a quasi-isomorphic one
  leads to a quasi-isomorphism between the whole constructions.

  For~$\sigma\in\Sigma$, let $T_\sigma$ be the torus corresponding
  to the sublattice of~$N$ spanned by the~$x_i$ not contained in~$\sigma$.
  Then $X_\sigma=Z_\sigma\times T_\sigma$
  for some contractible toric variety~$Z_\sigma$
  (\cf~\cite[Sec~2.3]{Fulton:93}). Therefore,
  the map~$C(X_\sigma)\to C(T_\sigma)$
  is a quasi-isomorphism of $\Ll$-coalgebras, natural in~$\sigma$.
  Replacing in~$E(\U)$ each~$C(X_\sigma)$ by~$C(T_\sigma)$,
  we obtain another complex~$\Sl\tildeotimes E'(\U)$,
  quasi-isomorphic to~$\Sl\tildeotimes E(\U)$ as $\Sl$-coalgebra.

  Let $\Sl_\sigma\subset\Sl$ be the subcoalgebra
  cogenerated by the~$x_i$ spanning~$\sigma$.
  We now want to substitute $\Sl_\sigma$
  for each term~$\Sl\tildeotimes C(T_\sigma)$
  appearing in~$\Sl\tildeotimes E(\U)$.
  First of all, the inclusion
  $$
    \Sl_\sigma\to\Sl\tildeotimes C(T_\sigma),
    \quad
    s\mapsto s\otimes1
  $$
  is a chain map because of our choice of representatives of the~$x_i$.
  Moreover, it is natural in~$\sigma$.
  To see that it induces an isomorphism in homology,
  filter $\Sl\tildeotimes C(T_\sigma)$ by $\Sl$-degree.
  Then the $E^1$~term of the associated spectral sequence
  is~$\Sl\tildeotimes\Ll^\sigma$, which is $\Sl_\sigma$
  tensored with a Koszul complex.
  Therefore, the $E^2$~term is~$\Sl_\sigma$, as desired.
  We thus obtain another double complex~$E''(\U)$,
  quasi-isomorphic as $\Sl$-coalgebra to the preceding one.

  An argument analogous to that for exactness of the Mayer--Vietoris
  spectral sequence (\cf~\cite[\S 15]{BottTu:82})
  shows that the canonical map~$E''(\U)\to\SRco\Sigma$
  given by the inclusion~$\Sl_\sigma\hookrightarrow\SRco\Sigma$
  for~$\sigma=\sigma_i$ a maximal cone (and the zero map otherwise)
  is a quasi-isomorphism because $\Sigma$ is regular.
  (To see this in cohomology, note that the sheaf on~$X/T$
  determined by~$X_\sigma/T\mapsto\Su_\sigma$ is flabby.)
  This finishes the proof.
\end{proof}

\begin{proposition}
  The coalgebra~$C(X_T)$ is quasi-isomorphic to~$\SRco\Sigma$.
\end{proposition}

\begin{proof}
  Let $Q\colon(N',\Sigma')\to(N,\Sigma)$ be the Cox construction~\cite{Cox:95}
  (which is analogous to the construction~$\mathcal Z_P$
  of Davis--Januszkiewicz~\cite{DavisJanuszkiewicz:91}).
  More precisely: Let $N'$ be the free $\Z$-module
  over the rays in~$\Sigma$. For each cone in~$\Sigma$
  define a cone in~$N'$ spanned by the same rays,
  now considered as elements of~$N'$. The collection of all these
  cones form a fan~$\Sigma'$, isomorphic to~$\Sigma$ as partially ordered set
  because $\Sigma$ is simplicial.
  The map~$Q$ sending each~$\rho\in N'$ to its generator in~$N$ is a morphism
  of fans, hence induces a morphism~$X'=X_{\Sigma'}\to X=X_\Sigma$
  of toric varieties. It will be convenient to assume that the latter
  is surjective. If this is not the case (\ie, if the rays of~$\Sigma$
  do not span $N$), we replace $N'$ by~$N'\oplus N$ (with the obvious map
  to~$N$) and consider $\Sigma$ as lying in this enlarged space.

  The map~$X'\to X$ is actually a principal fibration with fibre,
  the subtorus~$\bar T\subset T'=T_{N'}$ corresponding
  to the kernel~$\bar N$ of~$Q\colon N'\to N$.
  Hence, the map~$C(X'_{T'})\to C(X_T)$
  is a quasi-isomorphism of coalgebras.
  Moreover, $C(X'_{T'})$ and $\Sl'\tildeotimes C(X')$
  are quasi-isomorphic coalgebras as well
  by~\cite[Cor.~4.4 \& Thm.~4.7]{Franz:03}.
  Together with Proposition~\ref{subfan-positive-orthant-formal}
  this finishes the proof
  since $\SRco\Sigma$~and~$\SRco{\Sigma'}$ are isomorphic.
\end{proof}

\section{Proof of Theorem~\protect\ref{cohomology-algebra}}

Consider the complex
$$
  L:=\Ll\tildeotimes\Sl\tildeotimes C(X)
$$
with differential
$$
  (-1)^{\degree a}d(a\otimes s\otimes c)=a\otimes s\otimes d\,c
    +\sum_{i=1}^r\Bigl(
      a\wedge x_i\otimes\xi_i\cap s\otimes c
      -a\otimes\xi_i\cap s\otimes x_i\cdot c
    \Bigr).
$$
It is a $\Ll$-coalgebra
with $\Ll$ acting on the first factor and
the diagonal being as for ordinary tensor products
of coalgebras.
The projection
$$
  L\to C(X),
  \qquad
  a\otimes s\otimes c\mapsto\epsilon(s)\,a\cdot c
$$
is a quasi-isomorphism of $\Ll$-coalgebras.
Here $\epsilon\colon\Sl\to\Ring$ denotes the canonical augmentation
sending each cogenerator to~$0$.
($L$ is actually the Koszul dual of the $\Sl$-comodule~$\Sl\tildeotimes C(X)$,
and the map to~$C(X)$ is one of the canonical quasi-isomorphisms
relating the Koszul functors, \cf~\cite[Sec.~2.6]{Franz:03}.)
Now replace $\Sl\tildeotimes C(X)$
by the quasi-isomorphic $\Sl$-coalgebra~$\SRco\Sigma$.
This leads to the complex~$\Ll\tildeotimes\SRco\Sigma$,
which is quasi-isomorphic to~$L$ because this tensoring with~$\Ll$
preserves quasi-isomorphisms, \cf~\cite[Prop.~2.2]{Franz:03}.
The dual of~$\Ll\tildeotimes\SRco\Sigma$
was given in Theorem~\ref{cohomology-algebra}.

\section{Proof of Theorem~\protect\ref{cohomology-module}}

  Consider again the Cox construction~$X'\to X$
  (which we suppose to be surjective).
  Choose a basis of the kernel~$\bar N$ of the projection~$N'\to N$,
  complete it to basis~$(x'_i)$ of~$N'$ and choose representative loops
  in~$C(T')$ as in Section~\ref{toric}. This induces a basis of~$N$
  with representatives loops given by projection.
  We therefore obtain quasi-isomorphisms
  of algebras~$\Ll\to C(T)$ and~$\Llprime\to C(T')$ compatible
  with the maps $\Llprime\to\Ll$~and~$C(T')\to C(T)$.
  Denote the torus corresponding to~$\bar N$ by~$\bar T$
  and its homology by~$\bar\Ll$.

  By Moore's theorem (see~\cite[Thm.~7.27]{McCleary:01} for instance),
  there is an isomorphism (for the moment being, of $\Ring$-modules)
  \begin{equation}\label{tor-Moore}
    H(X)=\Tor^{C(\bar T)}(\Ring,C(X'))
  \end{equation}
  because $X'\to X$ is a principal $\bar T$-bundle.
  We replace $C(\bar T)$ by $\bar\Ll$
  and $C(X')$
  by the free~$\bar\Ll$-module~$L:=\Llprime\tildeotimes\Sl'\tildeotimes C(X')$
  as in the proof of Theorem~\ref{cohomology-algebra}.
  Then
  \begin{equation}\label{tor-Moore-replace}
    H(X)=\Tor^{\bar\Ll}(\Ring,L)
      =H(\Ring\otimesunder{\bar\Ll}
           \Llprime\tildeotimes\Sl'\tildeotimes C(X'))
    .
  \end{equation}
  Now we would like to substitute
  the Stanley--Reisner coalgebra~$\SRco{\Sigma'}$
  for~$\Sl'\tildeotimes C(X')$.
  But the basis~$(x'_i)$ of~$N'$ was chosen to be compatible with~$\bar N$
  and does not contain the rays of~$\Sigma'$ in general.
  In other words, the~$x_\rho$ do not act by loops in~$T'$
  around the respective factors of the decomposition of the torus
  corresponding to this basis.
  But this was crucial for the proof
  of Proposition~\ref{subfan-positive-orthant-formal}.

  By Proposition~\ref{change-representatives} below,
  the change of representatives extends to a quasi-isomorphism
  of $\Sl'$-comodules
  from~$\Sl'\tildeotimes C(X')$ to itself,
  but with different $\Ll'$-actions on~$C(X')$.

  Now that the representatives are right, we can replace
  $\Sl'\tildeotimes C(X_{\Sigma'})$ by~$\SRco{\Sigma'}$,
  which leads to the complex
  $$
    \Ring\otimesunder{\bar\Ll}\Llprime\tildeotimes\SRco{\Sigma'}
    =\Ll\tildeotimes\SRco\Sigma,
  $$
  the dual of which was given in Theorem~\ref{cohomology-module}.

  We finally turn to the $\Ll$-action.
  The maps
  $$
    \Ring\otimesunder{\bar\Ll}\Llprime\tildeotimes\SRco{\Sigma'}
    \to\Ring\otimesunder{\bar\Ll}\Llprime\tildeotimes\Sl'\tildeotimes C(X')
    \to\Ring\otimesunder{C(\bar T)}C(X')
    \to C(X)
  $$
  (which induce the isomorphisms
  \eqref{tor-Moore}~and~\eqref{tor-Moore-replace})
  are $\Llprime$-equivariant, where the left $\Llprime$-action
  on~$\Ring\otimesunder{C(\bar T)}C(X')$ comes from that of~$C(X')$.
  This is well-defined because all groups are commutative
  and because the shuffle map is associative and commutative.
  Since $\bar\Ll$ acts trivially, we have in fact equivariance
  with respect to~$\Ll$.

\begin{proposition}\label{change-representatives}
  Let $L$ be a $C(T)$-module, turned into a $\Ll$-module
  by a quasi-isomorphism~$\phi\colon\Ll\to C(T)$.
  Suppose that $\phi'\colon\Ll\to C(T)$ is another quasi-iso\-mor\-phism,
  and denote by~$L'$ the same $C(T)$-module,
  but now equipped with the $\Ll$-action defined by~$\phi'$.
  Then the identity on~$L$ lifts to a quasi-isomorphism
  of weak $\Ll$-modules~$L\to L'$ \cite[Sec~2.5]{Franz:03},
  \ie, there is a quasi-isomorphism
  of $\Sl$-comodules~$\Sl\tildeotimes L\to\Sl\tildeotimes L'$
  which is the identity on~$1\otimes L=1\otimes L'$.
\end{proposition}

\begin{proof}
  In addition to the chain map~$f_0=\id_L\colon L\to L'$,
  we have to find maps~$f_\alpha\colon L\to L'$, $0\ne\alpha\in\N^r$,
  of degrees~$2(\alpha_1+\cdots+\alpha_r)$ such that
  $$
    \forall\alpha\in\N^r,\;l\in L
    \qquad
    d f_\alpha(l)-f_\alpha(d l)=\sum_{\alpha_i>0}\Bigl(
      c'_i\cdot f_{\alpha-i}(l)-f_{\alpha-i}(c_i\cdot l)
    \Bigr),
  $$
  where $c_i=\phi(x_i)$, $c'_i=\phi'(x_i)$, and $\alpha-i$ denotes
  $\alpha$ with the $i$-th coordinate decreased by~$1$.

  Since $T$ is commutative, $C(T)$ is a divided power algebra.
  For~$\Ring=\Z$ this means that the quotient
  $$
    c^{[k]}=\frac1{k!}\,c^k
  $$
  is well-defined for any $c\in C(T)$ of even degree.
  (See~\cite[Ex.~V.6.5.4]{Brown:82} for a proof in the case
  of the simplicial construction of~$T$ as the classifying space of~$N$.
  The proof for arbitrary commutative groups is essentially the same.)
  For other~$R$ we use extension of scalars.

  By assumption, the cycles $c_i$~and~$c'_i$ are homologous.
  Choose $b_i\in C_2(T)$ with~$d b_i=c'_i-c_i$ and define for~$\alpha\in\N^r$
  $$
    f_\alpha(l)=b_1^{[\alpha_1]}\cdots b_r^{[\alpha_r]}\cdot l.
  $$
  Using the relation~$d b_i^{[k+1]}=b_i^{[k]}\cdot(c'_i-c_i)$,
  it is straightforward to check that these $f_\alpha$ define
  a morphism of weak $\Ll$-modules.
\end{proof}

\section{Proof of Proposition~\protect\ref{Borel-Moore-homology}}
\label{proof-Borel-Moore-homology}

Denote by~$\Lu_\sigma\subset\Lu$ the subalgebra generated by linear forms
vanishing on (the linear hull of)~$\sigma\in\Sigma$, and by~$\xi_\sigma$
the product of the~$\xi_\rho$ over all rays~$\rho$ of~$\sigma$.
Note that $\Lu_\sigma$ is the dual of~$\Ll^\sigma$ and
stable under the action of~$\Ll$.

The following result is similar
to~\cite[Lemma~7.10]{BuchstaberPanov:02}.

\begin{lemma}
  The inclusion
  $$
    A(\Sigma)=\bigoplus_{\sigma\in\Sigma}\Lu_\sigma\otimes\xi_\sigma
    \hookrightarrow\Lu\otimes\SR\Sigma
  $$
  is a quasi-isomorphism of $\Ll$-modules.
\end{lemma}

\begin{proof}
  $A(\Sigma)$ is stable under the action of~$\Ll$, and also
  under the differential: Take a~$\sigma\in\Sigma$.
  Since the differential~\eqref{complex-cohomology}
  is independent of the chosen basis of~$N$, we may assume that the
  first $k=\dim\sigma$ basis vectors are the rays of
  (the regular cone)~$\sigma$, hence annihilate all elements in~$\Lu_\sigma$.
  Then for~$\alpha\otimes\xi_\sigma\in A(\Sigma)$
  \begin{align*}
    d(\alpha\otimes \xi_\sigma)
    &= \sum_{i=1}^k x_i\cdot\alpha\otimes\xi_i\xi_\sigma
      +\sum_{i=k+1}^r x_i\cdot\alpha\otimes\xi_i\xi_\sigma \\
    &= \sum_{i=k+1}^r\sum_{(\sigma,\rho)\in\Sigma}
         \xi_i(x_\rho)\,x_i\cdot\alpha\otimes\xi_i\xi_\sigma\\
    &= \sum_{\sigma<(\sigma,\rho)}x_\rho\cdot\alpha\otimes\xi_\rho\xi_\sigma
  \end{align*}
  where $(\sigma,\rho)$ denotes the cone spanned by $\sigma$
  and the (possibly redundant) ray~$\rho$.
  (If $\sigma$~and~$\rho$ do not span a cone of~$\Sigma$,
  then $\xi_\rho\xi_\sigma=0$.)
  Since the~$\xi_i$ are dual to the~$x_i$, they vanish on~$\sigma$
  for~$i>k$. This implies that in the above sum it suffices to consider
  only those~$\rho$ not contained in~$\sigma$. Then $\sigma$ is a facet
  of~$\tau=(\sigma,\rho)$, and $\xi_\rho\xi_\sigma=\xi_\tau$.

  To see that the inclusion induces an isomorphism in homology,
  filter~$\Lu\otimes\SR\Sigma$ according to
  the number of distinct rays~$\rho_i$ appearing
  in an element~$\alpha\otimes\xi_{\rho_1}\cdots\xi_{\rho_k}$.
  The $E_0$~term of the associated spectral sequence is composed of
  terms
  $$
    \Lu\otimes\Ring[\xi_{\rho_1},\ldots,\xi_{\rho_k}]\cdot\xi_\sigma
    =\Lu_\sigma\otimes\bigwedge(\xi_{\rho_1},\ldots,\xi_{\rho_k})
      \otimes\Ring[\xi_{\rho_1},\ldots,\xi_{\rho_k}]\cdot\xi_\sigma.
  $$
  This is $\Lu_\sigma\otimes\xi_\sigma$, tensored with a Koszul complex.
  Hence the $E_2$~term is $A(\Sigma)$.
\end{proof}

We note from the preceding proof
that the differential on~$A(\Sigma)$ has the form
\begin{equation}\label{differential-A-Sigma}
  d(\alpha\otimes \xi_\sigma)
  = \sum_{\sigma\facet\tau}x_{\rho(\tau,\sigma)}\cdot\alpha\otimes\xi_\tau  
\end{equation}
where $\rho(\tau,\sigma)$ denotes the additional ray of~$\tau\supset\sigma$.

\smallbreak

The orientation chosen for each~$\sigma\in\Sigma$
corresponds to a generator~$\omega'_\sigma$
of~$\bigwedge^{\dim\sigma}\lin\sigma$.
Together with an orientation~$\omega_0\in\bigwedge^r N$,
this uniquely defines
a generator~$\omega_\sigma\in\bigwedge^{r-\dim\sigma}N/\lin\sigma$
by~$\omega'_\sigma\wedge\omega_\sigma=\omega_0$.
The image of~$\omega_\sigma$ in~$\Ll^\sigma_{r-\dim\sigma}$
is denoted by the same symbol.
We compare orientations as follows:
Let $\sigma$ be a facet of~$\tau$
and denote the additional ray of~$\tau$ by~$\rho$. We set
\begin{align*}
  \omega'_\tau &= \Or_{\sigma,\tau}x_\rho\wedge\omega'_\sigma,
\intertext{or, equivalently,}
  \omega_\sigma &= \Or_{\sigma,\tau}\omega_\tau\wedge x_\rho.
\end{align*}

The cap product between $\Lu$~and~$\Ll$ is defined by
$$
  \beta(\alpha\cap a)=(\beta\wedge\alpha)(a)
$$
for $a\in\Ll$~and~$\alpha$,~$\beta\in\Lu$.

\begin{proposition}
  The map
  $$
    \pi\colon A(\Sigma)\to C(\Sigma),
    \quad
    \alpha\otimes\xi_\sigma\mapsto\alpha\cap\omega_\sigma\in \Ll^\sigma
  $$
  is, up to a degree shift by~$2r$, an isomorphism of $\Ll$-modules.
\end{proposition}

(Recall that we grade cohomology negatively.)

\begin{proof}
  The map is clearly bijective, so we only have to show that it is
  compatible with the differentials.
  Take a~$\alpha\otimes\xi_\sigma\in A(\Sigma)$.
  Then
  \begin{subequations}\label{two-differentials}
  \begin{align}
    d\,\pi(\alpha\otimes\xi_\sigma)
    &= \sum_{\sigma\facet\tau}
         \Or_{\sigma\tau}(-1)^{\degree\alpha+\degree{\omega_\sigma}}
         \pr_\tau(\alpha\cap\omega_\sigma)\\
  \intertext{and, by~\eqref{differential-A-Sigma}}
    \pi(d(\alpha\otimes\xi_\sigma))
    &= \sum_{\sigma\facet\tau}x_{\rho(\tau,\sigma)}\cdot\alpha\cap\omega_\tau
  \end{align}
  \end{subequations}
  We consider each summand, \ie, each~$\sigma\facet\tau$, separately.
  Extend the~$x_{\rho'}$, $\rho'$~a ray of~$\tau$,
  to a basis of~$N$.
  We distinguish two cases: If $\alpha$ is a product
  of linear forms~$\xi_i\in\Lu_\sigma$ dual to this basis,
  then it may or may not contain~$\xi_\rho$, $\rho=\rho(\sigma,\tau)$.
  If not, then both parts of~\eqref{two-differentials} are zero.
  If~$\alpha=\beta\wedge\xi_\rho$, then
  $$
    (x_\rho\cdot\alpha)\cap\omega_\tau
    =(-1)^{\degree\beta}(\beta\wedge x_\rho\cdot\xi_\rho)\cap\omega_\tau
    =(-1)^{\degree\alpha-1}\beta\cap\omega_\tau
  $$
  and
  $$
    \Or_{\sigma\tau}\pr_\tau(\alpha\cap\omega_\sigma)
    =(-1)^{\degree{\omega_\tau}}
      (\beta\wedge\xi_\rho)\cap(x_\rho\wedge\omega_\tau)
    =(-1)^{\degree{\omega_\sigma}-1}\beta\cap\omega_\tau.
  $$
  Again we have equality in~\eqref{two-differentials}.
\end{proof}

This duality implies that the spectral sequence for~$\Hcld(X)$
degenerates on the $E^2$~level.

\section{Proof of Corollary~\protect\ref{Chow-injective}}

Looking at~$H(\Sigma)=H(C(\Sigma))$ as the $E^2$~term
of the spectral sequence for~$\Hcld(X)$
coming from the filtration of~$X$ by orbit dimension, the cycle map
from Chow groups to~$H(\Sigma)$ is the inclusion of the diagonal terms.
Since by Proposition~\ref{Borel-Moore-homology}
this spectral sequence degenerates on the $E^2$~level
for any coefficients, we see that tensored cycle map
$$
  A(X)\otimes\Z_m \to\Hcld(X;\Z)\otimes\Z_m \to \Hcld(X;\Z_m)
$$
is injective for all~$m$.
Since a morphism~$M\to M'$ of finitely generated $\Z$-modules
is a split injection if and only if
the induced map~$M\otimes\Z_m\to M'\otimes\Z_m$ is injective for all~$m$,
this proves the claim.

\bibliographystyle{plain}
\bibliography{abbrev,algtop}

\end{document}